\documentclass[12pt,leqno]{article}
\usepackage{amsmath,amssymb,amsthm,amscd,euscript,mathrsfs}
\input xy
\xyoption{all}

\newtheorem{lemma}{Lemma}

\newtheorem{corollary}{Corollary}
\newtheorem{definition}{Definition}
\newtheorem{proposition}{Proposition}
\newtheorem{theorem}{Theorem}

\usepackage{verbatim}
\usepackage{fancyhdr}

\def\be{\beta}

 \def\mN{\mathbb N}
 
 \def\mP{\mathbb P}

 \def\mZ{\mathbb Z}

\def\gM{\mathfrak m}

 \def\kN{\mathcal N}
 \def\kO{\mathcal O}

\def\kE{\mathcal E} 
\def\kF{\mathcal F}

\def\kI{\mathcal I}

\def\kL{\mathcal L}

\def\rH{\mathrm H}

\begin{document}

\title{ACM sheaves on a nonsingular quadric hypersurface  in $\mP^4_k$}
\author{Elena Drozd}
\date{}

\maketitle

\begin{abstract}
We prove that on a nonsingular quadric hypersurface $Q$ in $\mP^4_k$ even
CI liaison classes of ACM curves are in bijective correspondence with the stable
equivalence classes (up to shift in degree) of maximal Cohen-Macaulay graded 
modules over the coordinate ring $R$ of $Q$, which in turn, are in bijective 
correspondence with stable equivalence classes (up to shift in degree) of ACM 
sheaves on $Q$ . In the situation of a nonsingular quadric 
hypersurface $Q\in \mP^4_k$ work of Kn\"orrer \cite{knorrer} shows that there is a 
unique nonfree indecomposable MCM module over $R$. We also describe the
ACM sheaf corresponding to this MCM module and give its cohomology table and its 
Hilbert polynomial.
\end{abstract}


~\newline

In classifying algebraic space curves, the special class of Arithmetically
Cohen-Macaulay (ACM)curves plays a significant role. C.Peskine and 
L.Szpiro (\cite {PS}) showed that any ACM scheme of codimension two in 
$\mP^4$ is in the CI-liaison class of a complete intersection. In higher 
codimension this result is no longer true for CI-liaison.

In this work our attention is restricted to the ACM curves that lie on a 
nonsingular quadric hypersurface $Q$ in $\mP^4_k$

We prove that on a nonsingular quadric hypersurface $Q$ in $\mP^4_k$ even
CI liaison classes of ACM curves are in bijective correspondence with the stable equivalence 
classes (up to shift in degree) of Maximal Cohen-Macaulay graded modules over the coordinate ring
$R$ of $Q$, which in turn, are in bijective correspondence with stable equivalence classes (up to
shift in degree) of ACM sheaves on $Q$ (\ref{theorem:secmain}).
In the situation of a nonsingular quadric hypersurface $Q\in \mP^4_k$ work of Kn\"orrer \cite{knorrer}
shows that there is a unique nonfree indecomposable MCM module over $R$. In this work we describe the
ACM sheaf corresponding to this MCM module. We give its cohomology table and its 
Hilbert polynomial.  This information is very helpful in determining which ACM sheaf 
coresponds to a given curve.  Some examples are discussed in \cite{drozd2}.

Let $X$ be a complete (connected) arithmetically Gorenstein subscheme of $\mP_k^n$
of dimension $r\geq 2$ with a very ample line bundle $\kL$. Suppose that $H^1(X,\kL^n)=0$ for 
all $n$. (For example $X=$ a smooth hypersurface of $\mP_k^4$.
 Then $h^i\bigl(X,O_X(t)\bigr)=0 \text{ }$ $ \forall t, i=1,\dots,r-1$ and $X$ is subcanonical: 
$w_X\cong O_X(c)$ for some 
integer $c$ since $X$ is arithmetically Gorenstein.)
We start with describing relations between ACM curves on $Q$ and MCM modules on $R$. 

%

\begin{definition}
Let $Z$ be a subscheme of $X$. Then a resolution $$0\longrightarrow {\mathcal E} \longrightarrow {\mathcal L}
\longrightarrow {\mathcal I}_Z \longrightarrow 0$$ where ${\mathcal L}=\oplus O_X(-a_i)$ for some $a_i\in \mathbb Z$ and
$\mathcal E$ is a locally free sheaf on $X$ with 
$H^1_*({\mathcal E})=0$ is called an 
\emph{$\mathcal E$-type resolution} of ${\mathcal I}_Z$. For any sheaf $\kF$ we denote $H^1\bigl(X,\kF(t)\bigr)$ for all
$t\in \mZ$ by $H^1_*(\kF)$.
\end{definition}

\begin{lemma}\label{lemma:etype}
Let $Z\subset X$ be locally Cohen-Macaulay equidimensional and $\text{codim}_X Z=2$. Then there exists an 
$\mathcal E$-type resolution ${\mathcal I}_Z$.
\begin{proof}

Let $S(X)$ be the homogeneous coordinate ring of $X$ and $S(Z)$ be the homogeneous coordinate ring of $Z$,
let $I_Z$ be the ideal of $Z$ in $S(X)$, and let 
\begin{equation}\label{equation:seq1}
0 \longrightarrow E \longrightarrow L \xrightarrow{\alpha} I_Z \longrightarrow 0
\end{equation}
be
a presentation where $L$ is a free graded $S(X)$-module (i.e. $\oplus S(-a_i)$ ), whose generators are
sent by $\alpha$ onto the generators of the homogeneous ideal of $Z$ in $S(X)$. Let $E=\text{ker}\alpha$.
By sheafifying one obtains a resolution 
\begin{equation}\label{equation:seq2}
0 \longrightarrow {\mathcal E} \longrightarrow {\mathcal L} \longrightarrow 
\kI_Z \longrightarrow 0 \text{,}
\end{equation}
where $\mathcal L$ is a direct sum of
line bundles on $X$.
To show that $\mathcal E$ is locally free sheaf on $X$,
consider any closed point $y\in Z$. We know that 
$\text{proj}.\text{dim} {\mathcal O}_{Z,y}=\text{dim}X-\text{depth}{\mathcal O}_{Z,y}=2$, 
since ${\mathcal O}_{X,y}$ is
regular and $Z$ is locally Cohen-Macaulay. So, $\mathrm{proj.dim}\kI_Z=1$. This implies that ${\mathcal E}$ is locally free. 
To show that (\ref{equation:seq2}) is an $\kE$-type resolution of $\kI_Z$ we need to prove that 
$H^1_*(\kE)=0$. From (\ref{equation:seq2}) we get the following exact sequence:
\begin{equation}\label{equation:seq3}
0 \longrightarrow H^0_*(\kE) \longrightarrow H^0_*(\kL) \longrightarrow 
H^0_*(\kI_Z) \longrightarrow 
H^1_*(\kE) \longrightarrow H^1_*(\kL) \longrightarrow \dots .
\end{equation}
$H^1_*(\kL) =0$ since $\kL=\bigoplus_{i=1}^m\kO(-a_i)$. $H^0_*(\kL)=L$; $H^0_*(\kI_Z)=I_Z$. Therefore
sequence (\ref{equation:seq3}) is sequence (\ref{equation:seq1}). Thus $H^1_*(\kE)=0$ and sequence  
(\ref{equation:seq3}) is an $\kE$-type resolution of $\kI_Z$. This concludes the proof of the lemma \ref{lemma:etype}
\end{proof}
\end{lemma}

%
%

\begin{definition}
Two vector bundles $\mathcal F$ and $\mathcal G$ on $X$ are called \emph{stably equivalent} if 
$${\mathcal F}\oplus \bigoplus_{i=1}^t O_X(a_i)\cong {\mathcal G} \text{ } \oplus \bigoplus_{i=1}^s O_X(b_i)$$
for some $c,t,d,a_i,b_i\in \mathbb Z$.
\end{definition}

\begin{theorem}[\cite{rao} 1.9]\label{corollary:rao}
Let $X$ be a nonsingular qadric hypersurface in $\mP_k^4$. Then the even liaison classes of curves
$Z$ on $X$ are in bijective correspondence with the stable
equivalence classes modulo twists of vector bundles $\mathcal F$ on $X$ with the property that $H^1_*({\mathcal F})=0$.
The correspondence is given by $Z\longrightarrow \mathcal F$, where $\mathcal F$ is the kernel of an
$\mathcal E$-type resolution of ${\mathcal I}_Z$.
\end{theorem}

We note here that theorem \ref{corollary:rao} is the statement cited since the nonsingular quadric hypersurface
in $\mP^4_k$ is a complete connected Gorenstein scheme of dimension at least two. $\kO_X(1)$ is a 
very ample locally free sheaf with $H^i_*(\kO_X)=0$ for $i=1,2$. Finally, $\kL$-stable equivalence
of \cite{rao} 1.9 is stable equivalence up to shift in degree defined in this work.

\begin{lemma}\label{claim:ACM}
Let $C$ be a curve on a nonsingular hypersurface $Q$ in $\mP^4_k$. Let 
\begin{equation}\label{equation:seq2star}
0\longrightarrow \kF \longrightarrow\kL \longrightarrow\kI_C\longrightarrow 0
\end{equation}
be an $\kE$-type resolution of $\kI_C$. Then the curve $C$ is ACM if and only if $H^i_*(\kF)=0$ for $i=1,2$. 
\begin{proof}
From theorem \ref{corollary:rao} we know that $H^1_*(\kF)=0$. Recall that 
$C$ is ACM is equivalent to $H^1_*(\kI_C)=0$, where $\kI_C$ is the ideal sheaf of $C$.
Taking cohomology in (\ref{equation:seq2star}) we get an exact sequence:
$$H^1_*(\kF)\longrightarrow H^1_*(\kL)\longrightarrow H^1_*(\kI_C)\longrightarrow H^2_*(\kF) \longrightarrow H^2_*(\kL)\longrightarrow \dots ,$$
in which $H^1_*(\kI_C)=0$ since $C$ is ACM curve, and $H^2_*(\kL)=0$ since $\kL$ is a direct sum of locally free
sheaves. Therefore, if $C$ is ACM, $H^i_*(\kF)=0$ for $i=1,2$. Reversing this argument we get the converse of the
statement.
\end{proof}
\end{lemma}

\begin{definition}[\cite{eisenbud}]
Let $B$ be a local ring and $M$ be a finitely generated $B$-module. Then we say 
$M$ is a \emph{Maximal Cohen-Macaulay} module if $\text{\emph{depth}}M=\text{\emph{dim}}B$.
\end{definition}

\begin{definition}
On a nonsingular 3-fold $Q$ a locally free sheaf $\kF$ with the property that $H^i_*(\kF)=0$ for $i=1,2$ is 
called and \emph{ACM sheaf}.
\end{definition}

\begin{theorem}\label{theorem:secmain}
Let $Q$ be a nonsingular  hypersurface in $\mP_k^4$. Then the following classes are in bijective correspondence:
\begin{enumerate}
\item even CI-liaison classes of ACM curves
\item stable equivalence classes up to shift in degree of ACM sheaves on $Q$
\item stable equivalence classes up to shift in degree of Maximal Cohen-Macaulay graded modules over the ring $R=k[x_0,\dots,x_4]/I_Q$, 
where $I_Q$ is the ideal of $Q$.
\end{enumerate}
\begin{proof}
(1) $\leftrightarrow$ (2) is Theorem \ref{corollary:rao} plus Lemma \ref{claim:ACM}.

(1) to (3): \quad
For an ACM curve $C$ there exists the following resolution:
$$0\longrightarrow E \longrightarrow L \xrightarrow{\alpha} I_C \longrightarrow 0$$ 
where $I_C$ is the ideal of C, $L$ is a graded free $R$-module whose generators are sent by $\alpha$ onto the generators
of $I_C$ and $E=\text{Ker } \alpha$. Let $\gM$ be the irrelevant maximal ideal $R_+$ of $R$. We associate to the
curve $C$, the module $E$.

Then, in order to prove that $E$ is MCM we need to 
show that $\text{depth }E=\text{dim }R$. Now, $\text{dim }Q=3$,
therefore $\text{dim }R=4$; let $A_C$ be the coordinate ring of $C$. 
Then $\text{dim }A_C=2$. Consider the following short exact sequence:
$$0\longrightarrow I_C \longrightarrow R \longrightarrow A_C \longrightarrow 0$$
From it, taking local cohomology, we will calculate $\text{depth } I_C$:
$$\begin{array}{lll}
\text{depth }A_C=2 & \text{ since $C$ is ACM, therefore } & H^i_\gM(A_C)=0, i=0,1, H^2_\gM(A_C)\neq 0 \\
\text{depth }R=4, & \text{ therefore } & H^i_\gM(R)=0, i\leq 3 
\end{array}$$
Consider 
\begin{multline*}
0\longrightarrow H^0_\gM(I_C) \longrightarrow H^0_\gM(R) \longrightarrow H^0_\gM(A_C) \longrightarrow H^1_\gM(I_C) \longrightarrow \\
\longrightarrow H^1_\gM(R) \longrightarrow H^1_\gM(A_C) \longrightarrow H^2_\gM(I_C) \longrightarrow H^2_\gM(R) \longrightarrow H^2_\gM(A_C) \longrightarrow \\
\longrightarrow H^3_\gM(I_C) \longrightarrow H^3_\gM(R) \longrightarrow H^3_\gM(A_C) \longrightarrow \dots ~ .
\end{multline*}
Therefore $H^i_\gM(I_C)=0 \quad i=0,1,2, \quad H^3_\gM(I_C)\neq 0$, therefore $\text{depth }I_C=3$ by the cohomological interpretation of depth.

Similarly consider
$$0\longrightarrow E \longrightarrow L \longrightarrow I_C \longrightarrow 0$$
Taking local cohomology:
\begin{multline*}
0\longrightarrow H^0_\gM(E) \longrightarrow H^0_\gM(L) \longrightarrow H^0_\gM(I_C) \longrightarrow H^1_\gM(E) \longrightarrow \\
\longrightarrow H^1_\gM(L) \longrightarrow H^1_\gM(I_C) \longrightarrow H^2_\gM(E) \longrightarrow H^2_\gM(L) \longrightarrow H^2_\gM(I_C) \longrightarrow \\
\longrightarrow H^3_\gM(E) \longrightarrow H^3_\gM(L) \longrightarrow H^3_\gM(I_C) \longrightarrow H^4_\gM(E) \longrightarrow H^4_\gM(L) \longrightarrow \dots,
\end{multline*}

$$\begin{array}{lll}
\text{depth }L=4, & \text{ therefore } & H^i_\gM(L)=0, i\leq 3, H^4_\gM(E)\neq 0 \\
\text{depth }I_C=3, & \text{ therefore } & H^i_\gM(I_C)=0, i\leq 2, H^3_\gM(I_C)\neq 0
\end{array}$$

Therefore, $H^i_\gM(E)= 0, i=0,1,2,3, H^4_\gM(E)\neq 0$, therefore $\text{depth }E=4$, therefore $E$ is MCM as was to be shown.

(3) to (2): \quad
Let $E$ be a graded MCM $R$-module. Let ${\mathcal F}= \widetilde E$ be the associated sheaf over $Q$. To show that 
$H^i_\star({\mathcal F})=0 \quad\text{ for }i=1,2$ we need the following result:

\begin{theorem}[\cite{eisenbud}, p.693]
Let $R$ be a graded ring. If $M$ is a graded $R$-module, then there is a natural exact sequence (where $\gM$ is irrelevant 
maximal ideal $R_+$)
$$0\longrightarrow H^0_\gM(M) \longrightarrow M \longrightarrow H^0_*({\mathcal E}) \longrightarrow H^1_\gM(M) \longrightarrow 0$$
and $H^i_*({\mathcal E})\cong H^{i+1}_\gM(M)$ for $i>0$, where $\kE$ is $\widetilde M$.
\end{theorem} 

Applying this theorem for ${\mathcal F}= \widetilde E$ we obtain:
$$H^i_*{\mathcal F}\cong H^{i+1}_\gM(E) \quad i>0, \quad H^0_*({\mathcal F})=E,$$
which means that $H^1_*({\mathcal F})=H^2_\gM(E); H^2_*({\mathcal F})=H^3_\gM(E)$. Now $E$ is MCM, therefore $\text{depth }E=\text{dim }R=4$, 
therefore $H^i_\gM(E)=0, i<4$, therefore $H^2_\gM(E)=H^3_\gM(E)=0$, therefore $H^i_*({\mathcal F})=0, i=1,2$ as was to be shown. 

This concludes the proof of the theorem.
\end{proof}
\end{theorem}

\begin{corollary}\label{enciliaison}
Let $Q$ be a nonsingular  hypersurface in $\mP^4_k$. Let $C$ be an ACM curve on $Q$. 
If $\quad 0\to \kE\to \kL\to \kI_C\to 0$ is an $\kE$-type resolution of $\kI_C$, then there exists
a curve $C^\prime$ with an $\kE$-type resolution of $\kI_{C^\prime}$ of the form
$$0\to \kE\oplus \bigoplus_{i=1}^k\kO(a_i)\to \kL^\prime \to \kI_{C^\prime}\to 0 \text{, where }
\kL^\prime =\bigoplus_{j=1}^{k_1}\kO(b_j),$$
if and only if $C^\prime$ is CI-linked to $C$.

If $\quad 0\to \kL \to \kN \to \kI_C\to 0$ 
is an $\kN$-type resolution of $\kI_C$, then there exists a curve $C^\prime$
with an $\kN$-type resolution of $\kI_{C^\prime}$ of the form
$$0\to \kL^\prime \to \kN\oplus \bigoplus_{i=1}^k\kO(a_i)\to \kI_{C^\prime}\to 0 \text{, where }
\kL^\prime =\bigoplus_{j=1}^{k_1}\kO(b_j),$$
if and only if $C^\prime$ is CI-linked to $C$.
\end{corollary}

Now as equivalence of CI-biliaison classes of ACM curves and stable equivalence classes up to shift in degree 
of ACM sheaves is established, we describe ACM sheaves on a nonsingular quadric hypersurface $Q$ in $\mP^4_k$.

%
%

\begin{theorem}[\cite{BEH}]\label{theorem:t720}
Let $k=\overline{k}, \quad \text{\emph{char }}k\neq 2$. Let $Q$ be a quadratic form on a vector space $V$ over $k$ considered as an
element of $S_2(V^*)$, which is regular in the sense that $R={k[V^*]}/Q$ has only an isolated singularity. Then the Maximal 
Cohen-Macaulay modules over $R$ are such that there are always one or two nonfree indecomposable Maximal Cohen-Macaulay 
$R$-modules depending on whether $\text{\emph{dim }}V$ is odd or even; that they both have the same rank and are syzygies of
one another when there are two; and that writing $m=\frac{\text{\emph{dim }}V}{2}-1$, the rank of their direct sum is $2^m$.
\end{theorem}


\begin{corollary}\label{corollary:maxCM}
In the case of a nonsingular quadric hypersurface in $\mP^4_k$ 
there is only one indecomposable Maximal Cohen-Macaulay module over $R$
and its rank is $2$.
\end{corollary}

Now we describe this unique MCM module over $R$. 

\begin{proposition}\label{lineres}
Let $L$ be a line in $Q\subseteq \mP^4_k$. Then it has an $\kE$-type resolution 
of the form: 
\begin{equation}\label{equation:A}
0 \longrightarrow {\mathcal E} \longrightarrow {\mathcal O}^3_Q(-1) \longrightarrow {\mathcal I}_L \longrightarrow 0
\end{equation}
\begin{proof}
Let $I_L$ be an ideal of a line $L$ in $R$ where $R$ is the homogeneous coordinate ring of $Q$. Then there is an
 exact sequence:
$$
0 \longrightarrow E \longrightarrow R^3(-1) \xrightarrow{\be} I_L \longrightarrow 0
$$
where $E=\ker \be$. By sheafifying we obtain the following exact sequence:
\begin{equation}\label{equation:B}
0 \longrightarrow \kE \longrightarrow \kO_Q^3(-1) \longrightarrow \kI_L \longrightarrow 0
\end{equation}
Following the proof of \ref{lemma:etype} we get that (\ref{equation:B}) is an $\kE$-type resolution 
of $\kI_L$. 
By \ref{claim:ACM} this implies that $\kE$ is an
ACM sheaf on $Q$. Also note that $\mathrm{rank}\kE=2$.
\end{proof}
\end{proposition}

\begin{definition}\label{definition:e0}
We denote by $\kE_0$ the ACM rank 2 sheaf appearing in the $\kE$-type resolution of $\kI_L$. 
\end{definition}


\begin{theorem}\label{theorem:t725}
The maximal Cohen-Macaulay module $E$ corresponding to $\kE_0$  is the unique nonfree indecomposable 
MCM module on $R$, and $\kE_0$ is the unique
(up to shift) indecompsable ACM sheaf on $Q$ that is not isomorphic to a direct sum of line bundles.
\begin{proof}
Let us show that ${\mathcal E}_0$ is not isomorphic to a direct sum of two locally free sheaves.

Suppose ${\mathcal E}_0={\mathcal L_1 \oplus L_2}$, where ${\mathcal L}_i$ are locally free sheaves. By \ref{theorem:secmain}
this would imply
that a line on $Q$ is in the biliaison class of a complete intersection, 
which is not the case since $\text{deg }L$ is
odd while degree of a complete intersection on $Q$ is even and CI-liaison preserves parity.
Therefore ${\mathcal E}_0$ is not a direct sum of locally free sheaves.
Therefore ${\mathcal E}_0$ is an indecomposable ACM sheaf corresponding to an indecomposable MCM module.

By the theorem 
\ref{theorem:t720} there is only one isomorphism class of  nonfree indecomposable graded MCM modules over $R$.
Therefore $H^0_*({\mathcal E}_0)$ is the unique MCM given by that theorem.
\end{proof}
\end{theorem}

Theorem \ref{theorem:t725} and corollary \ref{corollary:maxCM} imply the follwing:

\begin{corollary}\label{corollary:fact}
Let $Q$ be a nonsingular quadric hypersurface in $\mP^4_k$. Let $\kE$ be an ACM sheaf on $Q$. Then
$$\kE=\bigoplus_{i=1}^{k_1}\kE_0(a_i)\oplus \bigoplus_{j=1}^{k_2}\kO_X(b_j)$$
for some $a_i, b_j \in \mZ$ and $k_1, k_2\in \mN$.
\end{corollary}

We now calculate dimensions of cohomology groups of $\kE_0$ and its Hilbert polynomial.

\begin{proposition}\label{e0cohom} ~

\begin{enumerate}
\item  The dimensions of cohomology groups of ${\mathcal E}_0$ are: 
$$
\begin{array}{c|c|c|c|c|c|c|c|c|c}
n & 0 & 1 & 2 & 3 & 4 & 5 & 6 & 7 & 8 \\
\hline
h^{0}\bigl(Q,{\mathcal E}_0(n)\bigr) & 0 & 0 & 4 & 16 & 40 & 80 & 140 & 224 & 336 
\end{array} ~ , $$
$$~ h^1_*(\kE_0)=h^2_*(\kE_0)=0 ~\text{ and } ~ h^0(\kE_0(n))=0 \text{ for } n<0;
$$
\item $\kE_0(n)$ is generated by global section for all $n\geq 2$;
\item $\kE_0^\vee=\kE_0(3)$.
\end{enumerate}
\begin{proof}
1) \quad Since $\kE_0$ is an ACM sheaf, we have $h^1_*(\kE_0)=h^2_*(\kE_0)=0$. 
Taking cohomology in the short exact sequence 
$$0 \longrightarrow {\mathcal I}_Q(n) \longrightarrow {\mathcal O}_{\mP_k^4}(n) \longrightarrow {\mathcal O}_Q(n) \longrightarrow 0$$
we arrive at:
$$0 \longrightarrow H^0\bigl(  {\mathcal I}_Q(n) \bigr)\longrightarrow H^0\bigl( {\mathcal O}_{\mP_k^4}(n) \bigr)\longrightarrow 
H^0\bigl( {\mathcal O}_Q(n) \bigr)\longrightarrow H^1\bigl({\mathcal I}_Q(n) \bigr)=0, $$
where $H^1\bigl(\kI_Q(n)\bigr)=0$ since $Q$ is ACM. Therefore dimensions are as follows:
$$\begin{array}{r}
n \\
 \\
<0 \\
0 \\
1 \\
2 \\
3 \\
4 \\
5 \\
6 \\
7 \\
8  
\end{array} \left ( \begin{array}{ccc}
h^0\bigl(\kI_Q(n)\bigr) & \quad & h^0\bigl(\kO_{\mP^4}(n)\bigr)  \\
 & & \\
0 & \quad & 0 \\
0 & \quad & 1 \\
0 & \quad & 5 \\ 
1 & \quad & 15 \\
5  & \quad & 35 \\
15 & \quad & 70 \\
35 & \quad & 126 \\
70 & \quad & 210 \\
126 & \quad & 330 \\
210 & \quad & 495
\end{array} \right ) \Longrightarrow \begin{array}{c}
 h^0\bigl(\kO_Q(n)\bigr) \\
 \\
0 \\
1 \\
5 \\ 
14 \\
30 \\
55 \\
91 \\
140 \\
204 \\
285  
\end{array}$$
Similarly, from the short exact sequence 
$$0 \longrightarrow {\mathcal I}_L \longrightarrow {\mathcal O}_Q \longrightarrow {\mathcal O}_L \longrightarrow 0$$
twisting and taking cohomology we get
$$0 \longrightarrow H^0\bigl({\mathcal I}_L(n)\bigr) \longrightarrow H^0\bigl({\mathcal O}_Q(n)\bigr) \longrightarrow 
H^0\bigl({\mathcal O}_L(n)\bigr) \longrightarrow H^1\bigl({\mathcal I}_L(n)\bigr)=0,$$
where $H^1\bigl(\kI_L(n)\bigr)=0$ since a line $L$ is an ACM curve. Then dimensions of cohomology groups are:
$$\begin{array}{r}
n \\
 \\
<0 \\
0 \\
1 \\
2 \\
3 \\
4 \\
5 \\
6 \\
7 \\
8  
\end{array} 
\begin{array}{c}
h^0\bigl(\kI_L(n)\bigr) \\
 \\
0 \\
0 \\
3 \\ 
11 \\
26 \\
50 \\
85 \\
133 \\
196 \\
296  
\end{array} \Longleftarrow
\left ( \begin{array}{ccc}
h^0\bigl(\kO_Q(n)\bigr) & \quad & h^0\bigl(\kO_L(n)\bigr)  \\
 & & \\
0 & \quad & 0 \\
1 & \quad & 1 \\
5 & \quad & 2 \\ 
14 & \quad & 3 \\
30  & \quad & 4 \\
55 & \quad & 5 \\
91 & \quad & 6 \\
140 & \quad & 7 \\
204 & \quad & 8 \\
285 & \quad & 9
\end{array} \right ) $$
And similarly from the short exact sequence 
\begin{equation}\label{qqq}
0 \longrightarrow {\mathcal E}_0 \longrightarrow {\mathcal O}_Q^3(-1) \longrightarrow {\mathcal I}_L \longrightarrow 0
\end{equation}
defining $\kE_0$, 
twisting and taking chomology we get 
$$0 \longrightarrow H^0\bigl({\mathcal E}_0 (n)\bigr)\longrightarrow H^0\bigl({\mathcal O}_Q^3(n-1) \bigr)\longrightarrow 
H^0\bigl({\mathcal I}_L (n)\bigr)\longrightarrow 0 $$
since $H^1\bigl({\mathcal E}_0 (n)\bigr)=0$ for the Maximal Cohen-Macaulay module ${\mathcal E}_0$. Then dimensions
of cohomology groups are:
$$\begin{array}{rcc}
n & \quad & h^0\bigl(\kE_0(n)\bigr) \\
 & & \\
<0 & \quad & 0 \\
0 & \quad & 0 \\
1 & \quad & 0 \\
2 & \quad & 4 \\
3 & \quad & 16 \\
4 & \quad & 40 \\
5 & \quad & 80 \\
6 & \quad & 140 \\
7 & \quad & 224 \\
8 & \quad & 336 \\
\end{array} 
\Longleftarrow
\left ( \begin{array}{ccc}
h^0\bigl(\kO^3_Q(n-1)\bigr) & \quad & h^0\bigl(\kI_L(n)\bigr) \\
 & & \\
0 & \quad & 0 \\
0 & \quad & 0 \\
3 & \quad & 3 \\ 
15 & \quad & 11 \\
42  & \quad & 26 \\
90 & \quad & 50 \\
165 & \quad & 85 \\
273 & \quad & 133 \\
420 & \quad & 196 \\
612 & \quad & 276
\end{array} \right ) $$

2) \quad Note that $\kE_0$ is 2-regular \cite[1.1.4]{migliore} and thus, by Castelnuovo-Mumford
regularity \cite[1.1.5 (3)]{migliore} $\kE_0(n)$ is generated by its global section for all $n\geq 2$.

3) \quad Since $\kE_0$ is a locally free sheaf of rank 2, $\kE_0^\vee=\kE_0(-c_1)$, where $c_1$ is the first
Chern class of $\kE_0$. We calculate $c_1$ from the exact sequence (\ref{qqq}):
$$c_1(\kE_0)+c_1(\kI_L)=c_1\big(\kO_Q^3(-1)\big)=-3.$$
$L$ is of codimension 2 on $Q$, thus $c_1(\kI_L)=0$, wherefrom $c_1(\kE_0)=-3$. 
\end{proof}
\end{proposition}


\begin{proposition}
The Hilbert polynomial of ${\mathcal E}_0$ is 
$$P_{{\mathcal E}_0}(n)=4\binom{n+1}{3}$$
\begin{proof}
Consider $$0\longrightarrow {\mathcal I}_Q(n) \longrightarrow {\mathcal O}_{\mP_k^4}(n) \longrightarrow {\mathcal O}_Q(n) \longrightarrow 0.$$ 
The Hilbert polynomial of 
${\mathcal O}_Q(n)$ is $$P_Q(n) = {n+4 \choose 4} - {n+2 \choose 4}.$$ From the short exact sequence 
$$0\longrightarrow {\mathcal I}_L \longrightarrow {\mathcal O}_Q \longrightarrow {\mathcal O}_L \longrightarrow 0$$ 
we get that the Hilbert polynomial of $\kI_L$ is
$$P_L(n) = {n+4 \choose 4} - {n+2 \choose 4}-n-1.$$
Thus from the sequence 
$$0 \longrightarrow {\mathcal E}_0 \longrightarrow {\mathcal O}_Q^3(-1) \longrightarrow {\mathcal I}_L \longrightarrow 0.$$
we get the Hilbert polynomial of ${\mathcal E}_0$:
$$P_{{\mathcal E}_0}(n) = 3 \left [ {n+3 \choose 4} - {n+1 \choose 4} \right ] - {n+4 \choose 4} + {n+2 \choose 4}+n+1=4{n+1 \choose 3}$$
\end{proof}
\end{proposition}

 \begin{lemma}\label{e0res}
 There is an  $\kE$-type resolution of $\kE_0$ of the form 
 $$ 
   0\longrightarrow \kE_0(-1)\longrightarrow \kO^4(-2)\longrightarrow \kE_0\longrightarrow 0
 $$ 
 \end{lemma} 
 \begin{proof}
The dimension of  $\rH^0(\kE_0(2))$ is 4 and $\kE_0(2)$ is generated by global sections
by proposition \ref{e0cohom}. Thus, there is a surjective map
 $$ 
   \kO^4\longrightarrow \kE_0(2)\longrightarrow 0.
 $$ 
Its kernel $\kE$ is a  rank two ACM sheaf. Thus, by \ref{corollary:fact} $\kE$ is equal to 
either $\kE_0(a)$ or $\kO(a)\oplus\kO(b)$ for some $a,b\in \mZ$. By \ref{e0cohom} we must have $\kE=\kE_0(1)$. 
Thus, the sequence $$0\longrightarrow \kE_0(1)\longrightarrow \kO^4\longrightarrow \kE_0(2)\longrightarrow 0$$ is exact
and it gives an $\kE$-type resolution of $\kE_0(2)$.
 \end{proof}

\end{document}